\documentclass[12pt]{article}

\usepackage{graphicx}
\usepackage{geometry}

\renewcommand{\baselinestretch}{1.0}

\usepackage{amssymb}
\usepackage{amsmath ,bm}

\newcommand{\bey}{\begin{eqnarray}}
\newcommand{\eey}{\end{eqnarray}}
\newcommand{\beq}{\begin{equation}}
\newcommand{\eeq}{\end{equation}}

\newcommand{\bs}{\boldsymbol}
\newcommand{\boldm}[1] {\mathversion{bold}#1\mathversion{normal}}

\newtheorem{thm}{Theorem}[section]
\newtheorem{lemma}{Lemma}[section]

\begin{document}

\vspace*{0in}

\begin{center}

{\large \bf Sparse maximum likelihood estimation for \\[0.1in] regression models}


\bigskip
\bigskip
Min Tsao
\\{\small Department of Mathematics \& Statistics, University of Victoria\\British Columbia, Canada V8W 2Y2 }

\end{center}

\bigskip

{
\noindent {\bf Abstract:} For regression model selection via maximum likelihood estimation, we adopt a vector representation of candidate models and study the likelihood ratio confidence region for the regression parameter vector of a full model. We show that when its confidence level increases with the sample size at a certain speed, with probability tending to one, the confidence region consists of vectors representing models containing all active variables, including the true parameter vector of the full model. Using this result, we examine the asymptotic composition of models of maximum likelihood and find the subset of such models that contain all active variables. We then devise a consistent model selection criterion which has a sparse maximum likelihood estimation interpretation and certain advantages over popular information criteria.

\bigskip


\vspace{0.1in}

}


\section{Introduction}
There are a variety of model selection criteria based on a diverse range of motivations. Among the most commonly used criteria, Akaike Information Criterion (AIC) (Akaike, 1973) is based on choosing the candidate model with the minimum Kullback–Leibler divergence to the unknown true model, and Bayesian Information Criterion (BIC) (Schwarz, 1978) is motivated by selecting the candidate model with the highest posterior probability. Although their motivations differ widely, AIC and BIC have a maximum likelihood connection in that they both select models of maximum likelihood of certain sizes. Hannan and Quinn criterion (Hannan and Quinn, 1979) and Bridge criterion (Ding, Tarokh, and Yang, 2018a) also select models of maximum likelihood of certain sizes.
The method of maximum likelihood, dating back to R. A. Fisher (Aldrich, 1997), has been a cornerstone of statistical science for more than a century. It has applications in many areas. However, in the important area of model selection, there have been no selection criteria derived or justified solely based on this method in spite of the aforementioned connection which suggests that it has a fundamental role to play. In this paper, we fill this gap by studying regression models of maximum likelihood using familiar tools from the maximum likelihood theory toolbox, in particular, the classic work of S. S. Wilks on the asymptotic distribution of the likelihood ratio (Wilks, 1938). We develop a maximum likelihood estimation based approach for regression model selection which complements and lends support to information criteria.

The lack of a maximum likelihood based criterion may be partly due to the fact that candidate models in a general model selection problem may come from different parametric families. Since likelihood of models from different families cannot be compared directly, we cannot use the method of maximum likelihood to develop a general-purpose model selection criterion, so we have relied on information-theoretic, Bayesian, and other ideas that can be used to compare models from different families. See Kadane and Lazar (2004), Claeskens and Hjort (2008), and Ding, Tarokh and Yang (2018b), among others, for comprehensive discussions on model selection methods and philosophies.
Regression model selection problem, however, is a special model selection problem where candidate models are all submodels of the same full model, so they are in the same parametric family. It may be viewed as an estimation problem for a regression parameter vector of a full model that is known to be sparse. The method of maximum likelihood can be adapted to handle such an estimation problem. Modern $\ell_1$ penalized regression (Hastie, Tibshirani and Wainright, 2015) generates a sequence of sparse estimates for the regression parameter vector automatically via its $\ell_1$ penalty. Such estimates are restricted maximum likelihood estimates, and the final selection of an estimate is usually done non-parametrically via cross-validation. Our strategy is to focus on models of maximum likelihood of different sizes, first study their composition and then develop a criterion to select one of them as our estimate of the unknown true model. 


In Section 2, we discuss assumptions needed for all regression models under consideration. In Section 3, we study the likelihood ratio confidence region for the parameter vector of a full regression model and obtain an asymptotic characterization of its contents. We use this asymptotic characterization to first study the composition of models of maximum likelihood of different sizes, and then devise a sparse maximum likelihood estimator for the parameter vector of the full model. We conclude with some remarks in Section 4.

\section{Model assumptions}      

For simplicity, we will use linear model (\ref{e1}) for illustration but results obtained are also valid for generalized linear models satisfying assumptions (A1), (A2), and (A3) below. Consider a linear model with $p$ predictor variables,
\begin{equation}
\mathbf{y}=\mathbf{X}\bm{\beta} +\bm{\varepsilon}, \label{e1}
\end{equation}
where $\mathbf{y} \in \mathbb{R}^n$ is the vector of responses, $\bm{\beta}\in  \mathbb{R}^{p+1}$ is the vector of regression parameters, $\mathbf{X}\in  \mathbb{R}^{n\times (p+1)}$ is the design matrix, and $\bm{\varepsilon}\sim N(\mathbf{0}, \sigma^2\mathbf{I}_{n\times n})$. 
We assume that $p$ is fixed and $n>p$ throughout this paper except when we discuss the case of $p\gg n$ in Section \ref{remarks}.

Denote by $\bs{\beta}^t$ the true regression parameter vector of a full model such as (\ref{e1}), and by $\hat{\bs{\beta}}$ its maximum likelihood estimator (MLE). Let $\ell({\bs{\beta}})$ be the log-likelihood of a $\bs{\beta}\in \mathbb{R}^{p+1}$ and $\lambda({\bs{\beta}})=-2\{ \ell({\bs{\beta}})- \ell(\hat{\bs{\beta}})\}$ be the likelihood ratio statistic. We make three asymptotic assumptions with respect to the sample size $n\rightarrow \infty$ for all models under consideration:
\begin{itemize}
\item[(A1)]  The MLE $\hat{\bs{\beta}}$ is consistent, i.e., $\hat{\bs{\beta}}\stackrel{p}{\longrightarrow}\bs{\beta}^t$.
\item[(A2)]  Wilks's Theorem holds under $H_0: \bs{\beta}^t = \bs{\beta}$, i.e., $\lambda({\bs{\beta}}) \stackrel{d}{\longrightarrow} \chi^2_{p+1}$ under $H_0$.
\end{itemize}

By (A2), for an $\alpha \in (0,1)$, a $100(1-\alpha)\%$ asymptotic confidence region for $\bs{\beta}^t$ is
\beq
{\cal C}_{1-\alpha}=\{ \bs{\beta} \in \mathbb{R}^{p+1}: \lambda({\bs{\beta}})\leq \chi^2_{1-\alpha, p+1} \},  \label{c_region}
\eeq
where $\chi^2_{1-\alpha, p+1}$ denotes the $(1-\alpha)$th quantile of the $\chi^2_{p+1}$ distribution. Since $\lambda({\bs{\beta}})\geq 0$ and $\lambda(\hat{\bs{\beta}})=0$, $\hat{\bs{\beta}}$ is in the center of ${\cal C}_{1-\alpha}$. Hence, $\max_{\bs{\beta}\in {\cal C}_{1-\alpha}} \| \bs{\beta} -\hat{\bs{\beta}}\|_2$ is a measure of the size of ${\cal C}_{1-\alpha}$. The third assumption is concerned about the size and the confidence level.

\begin{itemize}
\item[(A3)]  The confidence level of ${\cal C}_{1-\alpha}$ can go to one and the size of  ${\cal C}_{1-\alpha}$ can go to zero at the same time in the sense that there exists a monotone decreasing sequence $\{\alpha_n\}_{n=1}^\infty$ such that  $\alpha_n\rightarrow 0$ and $\max_{\bs{\beta}\in {\cal C}_{1-\alpha_n}} \| \bs{\beta} -\hat{\bs{\beta}}\|_2=o_p(1)$.
\end{itemize}

These assumptions are valid under commonly used conditions for the asymptotic normality of $\hat{\bs{\beta}}$ which can be found in the literature.
To illustrate this point with the linear model, two commonly used conditions for the asymptotic normality of $\hat{\bs{\beta}}$ for (\ref{e1}) is
\beq
\frac{1}{n} \mathbf{X}^T\mathbf{X} \rightarrow D \quad \mbox{and}
\quad \frac{1}{n} \max_{1\leq i \leq n}  \mbox{\boldm $x$}_i^T\mbox{\boldm $x$}_i \rightarrow 0,  \label{e_cond}
\eeq
where $\mbox{\boldm $x$}_i$ is the $i$th row of $\mathbf{X}$ and $D$ is a $(p+1)\times (p+1)$ positive definite matrix. 
Under (\ref{e_cond}), 
\beq
\sqrt{n}(\hat{\bm{\beta}} - \bm{\beta}^t)  \stackrel{d}{\longrightarrow} N(\mathbf{0}, \sigma^2 D^{-1}) \label{asynorm}
\eeq
which implies (A1). Under $H_0: \bs{\beta}^t = \bs{\beta}$, by (\ref{asynorm}) the Wald statistic $W(\bs{\beta})$ satisfies
\beq 
W(\bs{\beta})=(\bs{\beta}-\hat{\bs{\beta}})^T\mathbf{X}^T\mathbf{X}(\bs{\beta}-\hat{\bs{\beta}})/\hat{\sigma}^2 \stackrel{d}{\longrightarrow}  \chi^2_{p+1},  \label{e2}
 \eeq
 where $\hat{\sigma}^2$ is the mean squared error of the full model. A two-term Taylor series expansion of the likelihood ratio statistic $\lambda(\bs{\beta})$ shows that
 \beq
\lambda(\bs{\beta})
=W(\bs{\beta})+o_p(1). \label{e2.5}
\eeq
By Slusky's theorem, (\ref{e2}) and (\ref{e2.5}) imply that (A2) holds. Because ${\cal C}_{1-
\alpha_n}$ has no analytic expression, it is difficult to evaluate its size directly. We now show that (A3) holds through the asymptotic equivalence between $\lambda(\bs{\beta})$ and $W(\bs{\beta})$. For a fixed $\gamma\in (0, 1)$, let 
\beq 
\alpha_{n}=1-P(\chi^2_{p+1} \leq n^\gamma). \label{e3}
\eeq
Then, $\alpha_n$ is monotone decreasing and converges to zero. Also, $\chi^2_{1-\alpha_n, p+1}=n^\gamma$. It follows from this and (\ref{e2.5}) that
$\lambda(\bs{\beta}) \leq \chi^2_{1-\alpha_n, p+1} = n^\gamma$ is equivalent to
\[ 
W(\bs{\beta})+o_p(1) \leq n^\gamma.
\]
Dividing both sides of the above inequality by $n$ shows that when $\lambda(\bs{\beta}) \leq  \chi^2_{1-\alpha_n, p+1}$, the corresponding $W(\bs{\beta})$ satisfies $n^{-1}W(\bs{\beta})=n^{-1}(\bs{\beta}-\hat{\bs{\beta}})^T\mathbf{X}^T\mathbf{X}(\bs{\beta}-\hat{\bs{\beta}})/\hat{\sigma}^2=O_p(n^{\gamma-1})$. This and the first condition in (\ref{e_cond}) imply (A3) holds for the sequence of $\alpha_n$ defined in (\ref{e3}). Different $\gamma$ values will generate different sequences of $\alpha_n$ in (\ref{e3}), and we will see later that they all lead to consistent model selection criteria. Similar phenomena have appeared in information criteria where there may be a range of penalty terms that all guarantee consistency of the criteria; see, for example, Rao and Wu (1989). Incidentally, $\alpha_n$ also plays the role of the penalty term in information criteria in that it controls the sparsity level of the selected model.  

For generalized linear models, (A1), (A2), and (A3) also hold under commonly used conditions for the asymptotic normality of $\hat{\bs{\beta}}$. See 
Haberman (1977), Gourieroux and Monfort (1981), and Fahrmeir and Kaufmann (1985) for such conditions.

\section{Main results}      

We first discuss a $(p+1)$-vector representation of a candidate model as our subsequent analyses will be based on vector representations of models. Suppose the intercept term $\beta_0$ is always included in a model. With $p$ predictor variables $\{x_1, x_2, \dots, x_p\}$, there are $2^p$ candidate models ${\cal M}=\{M_j\}_{j=1}^{2^p}$ where $M_j$ is defined as the $j$th subset of $\{x_1, x_2, \dots, x_p\}$. Each $\bs{\beta} \in \mathbb{R}^{p+1}$ represents a candidate model through the locations of zeros in its last $p$ elements which we will refer to as the sparsity structure of $\bs{\beta}$. Vectors with the same sparsity structure represent the same model, e.g.,  $\bs{\beta}_1=(1, 2, 0, 3, 0, \dots, 0)^T$ and $\bs{\beta}_2=(3, 1, 0, 2, 0, \dots, 0)^T$  both represent model $\{x_1, x_3\}$. Denote by $p^*$ the number of active variables (variables with non-zero parameters in the full model), and by $M^t_j$ the true model consisting of the $p^*$ active variables. Excluding $\beta_0$, non-zero elements of ${\bs{\beta}}^t$ represent variables in  $M^t_j$, and zeros of ${\bs{\beta}}^t$ represent inactive variables. 
Among vectors representing $M_j$, we denote the one with the highest likelihood by $\hat{\bs{\beta}}_j$ and call it the MLE for $M_j$. When $M_j$ is not the full model, $\hat{\bs{\beta}}_j$ contains zeros representing variables not in $M_j$. Let $\hat{\bs{\beta}}^t_j$ be the MLE for $M^t_j$. Then, $\hat{\bs{\beta}}^t_j$ has the same sparsity structure as ${\bs{\beta}}^t$ but a higher likelihood than ${\bs{\beta}}^t$. 

\subsection{Sparse maximum likelihood estimators} 
The MLE for the full model $\hat{\bs{\beta}}$ is a continuous random vector whose elements are, with probability one, all non-zero. When ${\bs{\beta}}^t$ is known to be sparse and its sparsity structure is of interest, $\hat{\bs{\beta}}$ is no longer a suitable estimator for ${\bs{\beta}}^t$ as it does not provide information about the sparsity structure of ${\bs{\beta}}^t$. The MLE for the true model $\hat{\bs{\beta}}^t_j$ is the ideal sparse estimator for ${\bs{\beta}}^t$ with the same sparsity structure as ${\bs{\beta}}^t$, but $\hat{\bs{\beta}}^t_j$ is unknown. We now explore known sparse maximum likelihood estimators. For a fixed $j\in\{0,1,\dots, p\}$, there are $p \choose j$ models with $j$ variables. Let $M_j^*$ be the model of maximum likelihood among these $p \choose j$ models in the sense that its MLE, denoted by $\hat{\bs{\beta}}^*_j$, has the highest likelihood among MLE's for these $p \choose j$ models. We call  $M_j^*$ and $\hat{\bs{\beta}}^*_j$ the $j$th sparse maximum likelihood estimates for $M^t_j$ and $\bs{\beta}^t$, respectively. We also call the collection of $(p+1)$ models ${\cal M}_{ml}=\{M^*_j\}_{j=0}^p$ the maximum likelihood set. Since there is only one empty model and only one full model, $M^*_0$ is the empty model and $M^*_p$ is the full model. Intuitively, we expect that models with $p^*$ or more variables in ${\cal M}_{ml}$ to be correct estimates for $M^t_j$ in the sense that they each contain all variables in $M^t_j$. We now show that this is asymptotically true. We need the following lemmas. Proofs of the lemmas are in the Appendix.

\begin{lemma} \label{lem1}
Assume (A2) holds. For a monotone decreasing sequence $\alpha_n$  that converges to 0, we have
\beq 
\lim_{n\rightarrow \infty} P( \bs{\beta}^t \in  {\cal C}_{1-\alpha_n} ) =1. \label{lemma1}
\eeq
\end{lemma}

Suppose there is at least one active variable ($p^*\geq 1$); that is, among the last $p$ elements of $\bs{\beta}^t$, there is one or more elements that are non-zero. Let $\beta^s$ be the smallest one (in absolute value) among such non-zero elements of $\bs{\beta}^t$. Denote by ${\cal B}$ the collection of $\bs{\beta} \in \mathbb{R}^{p+1}$ representing models with one or more active variables missing; that is, all $\bs{\beta}\in {\cal B}$ have one or more zeros in their elements for the $p^*$ active variables. Then, 
\[
\inf_{\mathbf{{\cal B}}}\| \bs{\beta} -\bs{\beta}^t\|_2= |\beta^s|>0.
\]
For a fixed $\delta$ such that $0<\delta < |\beta^s|$, define a neighbourhood of $\bs{\beta}^t$ as follows,
\beq
{\cal N}(\bs{\beta}^t, \delta) = \{ \bs{\beta} \in \mathbb{R}^{p+1}: \| \bs{\beta} -\bs{\beta}^t\|_2 \leq \delta\}.
\eeq
It is clear that ${\cal B} \cap {\cal N}(\bs{\beta}^t, \delta)=\phi$ which shows that ${\cal N}(\bs{\beta}^t, \delta)$ consists of $\bs{\beta}$ vectors representing models each containing all active variables. Lemma \ref{lem2} below gives the asymptotic relationship between ${\cal N}(\bs{\beta}^t, \delta)$ and ${\cal C}_{1-\alpha_n}$.

\begin{lemma} \label{lem2}
Assume that (A1) and (A3) hold, and that there is at least one active variable ($p^*\geq 1$). For the monotone decreasing sequence $\alpha_n$ in (A3), we have
\beq 
\lim_{n\rightarrow \infty} P\{ {\cal C}_{1-\alpha_n} \subset {\cal N}(\bs{\beta}^t, \delta) \} = 1. \label{lemma2}
\eeq
\end{lemma}

Combining Lemma \ref{lem1} and Lemma \ref{lem2}, we to obtain
\beq
\lim_{n\rightarrow \infty}  P({\bs{\beta}}^t \in  {\cal C}_{1-\alpha_n} \subset {\cal N}(\bs{\beta}^t, \delta) ) = 1, \label{rt5}
\eeq
which implies that, with probability tending to one, the likelihood ratio confidence region ${\cal C}_{1-\alpha_n}$ contains only vectors representing models with all active variables including ${\bs{\beta}}^t$. Note that (i) $\hat{\bs{\beta}}^t_j$ has a higher likelihood than ${\bs{\beta}}^t$, so $\hat{\bs{\beta}}^t_j$ is in ${\cal C}_{1-\alpha_n}$ if ${\bs{\beta}}^t$ is, (ii) since $\hat{\bs{\beta}}^t_j$ is the MLE of a model of $p^*$ variables, its likelihood is not higher than that of $\hat{\bs{\beta}}^*_{p^*}$  and thus $\hat{\bs{\beta}}^*_{p^*}$  is in ${\cal C}_{1-\alpha_n}$ if $\hat{\bs{\beta}}^t_j$ is, and (iii) for any $k> p^*$, $\hat{\bs{\beta}}^*_k$ has higher likelihood than $\hat{\bs{\beta}}^*_{p^*}$, thus $\hat{\bs{\beta}}^*_k$ is in ${\cal C}_{1-\alpha_n}$ if $\hat{\bs{\beta}}^*_{p^*}$ is. It follows from these and (\ref{rt5}) that, with probability tending to one, 
\beq
\{ \hat{\bs{\beta}}^*_{p^*}, \hat{\bs{\beta}}^*_{p^*+1}, \dots, \hat{\bs{\beta}}^*_{p}\}\subset {\cal C}_{1-\alpha_n} \label{rt10}
\eeq
and the corresponding models $\{ M^*_{p^*}, M^*_{p^*+1}, \dots, M^*_{p}\}$ each contains all active variables. Theorem \ref{thm1} below summarizes the above discussion.

\begin{thm} \label{thm1}
 Under (A1), (A2), and (A3), with probability tending to 1, a model with $p^*$ or more variables in the maximum likelihood set contains all active variables.
 \end{thm}
 
Theorem \ref{thm1} shows that asymptotically (a) $M^*_{p^*}, M^*_{p^*+1}, \dots, M^*_p$ are all correct estimates for $M^t_j$, and (b) $M^*_{p^*}=M^t_j$ because they are both the model with $p^*$ variables that contains all $p^*$ active variables. Point (b) implies that, for consistency in model selection, it suffices to consider only models in ${\cal M}_{ml}$ as the true model $M^t_j$ is in ${\cal M}_{ml}$. This reduces the dimension of the model space from $2^p$ to $(p+1)$. The above proof of Theorem \ref{thm1} implicitly assumed $p^*\geq 1$ as Lemma \ref{lem2} needs this condition, but the theorem holds trivially for $p^*=0$. 

\subsection{Constrained minimum criterion} 
Since $M^*_{p^*}=M^t_j$ with probability tending to one, $M^*_{p^*}$ is an ideal estimator for $M^t_j$ but unfortunately we do not know which model in ${\cal M}_{ml}$ is $M^*_{p^*}$ because $p^*$ is unknown. We now turn to sparse estimation of $\bs{\beta}^t$. The ideal estimator is $\hat{\bs{\beta}}^*_{p^*}$ which is unknown.  However, equation (\ref{rt10}) suggests that among sparse maximum likelihood estimators that are in the confidence region ${\cal C}_{1-\alpha_n}$, $\hat{\bs{\beta}}^*_{p^*}$ is the most sparse one. We thus use the most sparse maximum likelihood estimator in ${\cal C}_{1-\alpha_n}$ in place of the unknown $\hat{\bs{\beta}}^*_{p^*}$ to estimate $\bs{\beta}^t$.  Specifically, using the $\ell_0$ norm that counts the number of non-zero elements of a vector, denote by $\hat{\bs{\beta}}_n$ be the solution of the following minimization problem
\beq
\underset{\bm{{\cal M}_{ml}}}{\text{minimize}}       \|\hat{\bm{\beta}}_j^*\|_0
\mbox{\hspace{0.2in} subject to \hspace{0.01in} } \hat{\bm{\beta}}_j^* \in {\cal C}_{1-\alpha_n},  \label{def}
\eeq
and denote by $\hat{M}_n$ the model represented by $\hat{\bs{\beta}}_n$, we choose $\hat{\bs{\beta}}_n$ as the estimator for $\bs{\beta}^t$ and thus select model $\hat{M}_n$. We refer to this as the constrained minimum criterion (CMC) for regression model selection. We call $\hat{\bs{\beta}}_n$ the CMC estimator and $\hat{M}_n$ the CMC selection. Since the $(p+1)$ models in ${\cal M}_{ml}$ differ in their $\ell_0$ norm values and $\hat{\bm{\beta}}_p^* = \hat{\bm{\beta}} \in {\cal C}_{1-\alpha_n}$ regardless the confidence level $(1-\alpha_n)$, the solution $\hat{\bs{\beta}}_n$ for (\ref{def}) exists and is unique.

Tsao (2024) gave an earlier version of CMC with a fixed $\alpha$ level for selecting linear and generalized linear models,
and showed that the probability that the corresponding CMC selection is the true model is bounded below by $(1-\alpha)$ as $n$ goes to infinity. However, selection consistency cannot be achieved if $\alpha$ is fixed. In the present paper,  we allow $\alpha$ to depend on $n$ as described in assumption (A3). We now show that the resulting CMC estimator and CMC selection are both consistent. By the triangle inequality
\beq
\| \hat{\bs{\beta}}_n -\bs{\beta}^t\|_2 \leq \| \hat{\bs{\beta}}_n -\hat{\bs{\beta}}\|_2 + \| \hat{\bs{\beta}}-\bs{\beta}^t \|_2. \label{t.ineq}
\eeq
Under (A1) and (A3),  $\| \hat{\bs{\beta}}-\bs{\beta}^t \|_2=o_p(1)$ and $\max_{\bs{\beta}\in {\cal C}_{1-\alpha_n}} \| \bs{\beta} -\hat{\bs{\beta}}\|_2=o_p(1)$.  Since $\hat{\bs{\beta}}_n\in {\cal C}_{1-\alpha_n}$, these and (\ref{t.ineq}) imply that $\| \hat{\bs{\beta}}_n -\bs{\beta}^t\|_2=o_p(1)$. Thus, $\hat{\bs{\beta}}_n$ is a consistent estimator for ${\bs{\beta}}^t$.
To see that $\hat{M}_n$ is also consistent, we examine the following two cases. 

Case I: $p^*=0$.  For this case, $\bs{\beta}^t=(\beta_0, 0, \dots, 0)^T$ and $\hat{\bs{\beta}}_j^t=\hat{\bs{\beta}}_0^*=(\hat{\beta}_0, 0, \dots, 0)^T$. Since $\hat{\bs{\beta}}_j^t$ has a higher likelihood than $\bs{\beta}^t$, $P(\hat{\bs{\beta}}^t_j \in {\cal C}_{1-\alpha_n}) > P({\bs{\beta}}^t \in {\cal C}_{1-\alpha_n})$. By Lemma \ref{lem1}, $P({\bs{\beta}}^t \in {\cal C}_{1-\alpha_n})\rightarrow 1 $ which implies $P(\hat{\bs{\beta}}_0^* \in {\cal C}_{1-\alpha_n})=P(\hat{\bs{\beta}}^t_j \in {\cal C}_{1-\alpha_n}) \rightarrow 1$. Also, $\|\hat{\bs{\beta}}_0^*\|_0<\|\hat{\bs{\beta}}_j^*\|_0$ for $j\geq 1$, thus $\hat{\bs{\beta}}_n=\hat{\bs{\beta}}_0^*$ if $\hat{\bs{\beta}}^*_0 \in {\cal C}_{1-\alpha_n}$. It follows that $P(\hat{\bs{\beta}}_n=\hat{\bs{\beta}}_0^*)\rightarrow 1$ when $p^*=0$. 

Case II: $p^*\geq 1$. For this case Lemma \ref{lem2} implies that, with probability tending to one, none of $\hat{\bs{\beta}}^*_0, \dots,  \hat{\bs{\beta}}^*_{p^*-1}$ will be in ${\cal C}_{1-\alpha_n}$ as they represent models with fewer than $p^*$ variables which cannot contain all $p^*$ active variables. On the other hand, by (\ref{rt10}) we have $P(\{\hat{\bs{\beta}}^*_{p^*}, \hat{\bs{\beta}}^*_{p^*+1}, \dots, \hat{\bs{\beta}}^*_p\}\subset {\cal C}_{1-\alpha_n}) \rightarrow 1$. Also, $\|\hat{\bs{\beta}}_{p^*}^*\|_0<\|\hat{\bs{\beta}}_j^*\|_0$ for $j> p^*$. It follows that $\hat{\bs{\beta}}_n=\hat{\bs{\beta}}^*_{p^*}$, which implies $\hat{\bs{\beta}}_n=\hat{\bs{\beta}}^t_j$ and $\hat{M}_n={M}^t_j$, with probability tending to one. 

 Theorem \ref{thm2} below follows immediately from the above discussion.
 
\begin{thm} \label{thm2}

Under assumptions (A1),  (A2), and  (A3), we have (i) the CMC is estimation consistent in that
\beq \hat{\bm{\beta}}_n \stackrel{p}{\longrightarrow} \bm{\beta}^t, \label{convg1} \eeq
and (ii) the CMC is selection consistent in that
\beq P(\hat{M}_n= M_j^t) \rightarrow 1. \label{convg2} 
\eeq
\end{thm}

Assumption (A3) is not specific about the choice of $\alpha_n$. As long as it meets the conditions in (A3), the resulting CMC will be consistent. In finite sample applications, however, $n$ is fixed and the choice of $\alpha_n$ is crucial for serving different priorities of model selection. For a fixed $n$, when we increase the value of $\alpha_n$ from zero to one, the false active rate of the CMC selection goes up from zero to one and the false inactive rate goes down from one to zero. If $n$ is large relative to $p$, we recommend a small $\alpha_n$ (e.g., $0.01$, $0.05$, or $0.10$) so that both rates will be low due to the selection consistency of CMC. If $n$ is small or moderately large, it may not be possible to keep both rates low at the same time through the choice of $\alpha_n$, so we need to have priorities in order to choose an $\alpha_n$. If the priority is a low false inactive rate, based on our simulation study with multiple $\alpha_n$ levels, we recommend a large $\alpha_n$ (e.g., $0.8$ or $0.9$). The resulting CMC selection tends to be a large correct model in the maximum likelihood set containing all active variables but also some inactive variables. If the priority is a low false active rate or a low false discovery rate, we recommend a small $\alpha_n$ (e.g., 0.01, 0.05, or 0.10) so that the CMC selection is a small model from the maximum likelihood set which tends to contain few inactive variables but it may miss out some active variables. In simulations, we also observed that when the sample size is not large, CMC at level $0.5$ often gives the most balanced selection in that the two rates of the selected model are most comparable.

\subsection{Comparison with information criteria}  
Although they all select models from the maximum likelihood set,
information criteria use penalization to the log-likelihood to select models, whereas CMC uses the asymptotic $\chi^2$ calibration of the likelihood ratio which gives the CMC selection a familiar frequentist interpretation; the CMC selection at level $\alpha_n$ is the most sparse model in ${\cal M}_{ml}$ not rejected by the likelihood ratio test at level $\alpha_n$. Models selected by information criteria do not have a unified interpretation as their interpretations depend on the underlying motivations of the criteria which are different.

Simulation results show that, among AIC, BIC, and CMC at various $\alpha_n$ levels, CMC at level 0.5 is often the most accurate in terms of the overall misclassification rate (sum of false active rate and false inactive rate) of the selected model when the sample size is not very large.  To improve the selection accuracy of AIC and BIC in small sample situations, finite sample corrections to their penalty terms have been proposed by several authors, e.g., Hurvich and Tsai (1989), Broersen (2000) and Sclove (1987). CMC does not need such corrections. The $\chi^2$ calibration of the likelihood ratio offers a wide range of CMC criteria defined by a well-understood parameter $\alpha_n$ from the likelihood ratio test. We can simply choose CMC with different $\alpha_n$ levels to handle different situations and serve different priorities.

For selecting Gaussian linear models, BIC and CMC are both consistent, so we expect their probabilities of selecting the true model to be close to one when the sample size $n$ is large. When we set $n$ to very large values in our simulation study, the probabilities of CMC at levels $0.1, 0.05$ and $0.01$ are indeed very close to one, but that of the BIC is not as close to one because of a small but persistent false active rate. We observed that the likelihood ratios of models selected by BIC tend to be smaller than that of the CMC selections at the above three $\alpha_n$ levels. The small but persistent false active rate of the BIC is perhaps due to its tendency to pick up models with smaller likelihood ratios which are larger models that are more likely to contain inactive variables even when $n$ is very large.


Numerical evidence supporting the above observations may be found in Tsao (2024) and thus not repeated here. A comprehensive numerical study containing more supporting evidence as well as examples of CMC in high-dimensional situations will appear elsewhere.

\section{Concluding remarks} \label{remarks}

The asymptotic composition of regression models of maximum likelihood revealed in Theorem \ref{thm1} adds a new result particularly useful for model selection to the maximum likelihood theory. It supports information criteria and CMC for their selecting models of maximum likelihood. It also supports post-selection asymptotic inference for parameters of variables in their selected models conditioning on the models being sufficiently large (correct) models of maximum likelihood. Under this condition, the MLE of a selected model is consistent and asymptotically normally distributed, so we can construct confidence intervals for the parameters and test hypotheses about them as usual.

Theorem \ref{thm1} gives only a partial characterization of the maximum likelihood set that covers models of $p^*$ or more variables. Intuitively, as the sample size goes to infinity, models with 1 to $(p^*-1)$ variables may contain only active variables, but a proof of this seems to be beyond the reach of tools that we have used in this paper.
We used the likelihood ratio confidence region to define CMC in (\ref{def}). One may consider other types of confidence regions such as Wald confidence region. For Gaussian linear models, an $F$-statistic based confidence region can also be used (Tsao, 2021). Nevertheless, consistency of CMC based on other types of confidence regions is in general more difficult to establish due to the loss of the argument that ${\bs{\beta}}^t \in {\cal C}_{1-\alpha_n}$ implies $\hat{\bs{\beta}}^t_j \in {\cal C}_{1-\alpha_n}$ which is key to the proof of Theorem \ref{thm1} and only valid when the likelihood ratio confidence region is used.
There is also numerical evidence that the CMC criterion based on the likelihood ratio confidence region is more accurate.

Finally, we briefly consider the high-dimensional scenario of $p\gg n$. Suppose we may assume that $p^*\leq k$ where $k$ is given and $n>k$. Then, we only need to consider models of $k$ or fewer variables. Denote by ${\cal M}^k$ the set of all such models. There is now a fast algorithm that can handle best subset selection from ${\cal M}^k$ for $n$ in the hundreds and $p$ in the thousands (Bertsimas, King and Mazumder, 2016). MLE's of models in ${\cal M}^k$ are well-defined because $n>k$. Let ${\cal M}^k_{ml}=\{M^*_j\}_{j=0}^k$ be the set of the first $(k+1)$ sparse maximum likelihood estimators. Suppose $p$ is fixed and (A1), (A2), and (A3) hold as $n$ goes to infinity. Then, Theorem \ref{thm1} implies that, with probability tending to one, models with $p^*$ to $k$ variables in ${\cal M}^k_{ml}$ contain all active variables. This provides an asymptotic justification for selecting from ${\cal M}^k_{ml}$. This justification is valid even though CMC cannot be applied to select a model from ${\cal M}^k_{ml}$ because the likelihood ratios of candidate models are undefined for $p>n$. To be able to use CMC, we may remove from consideration variables that do not appear in models in ${\cal M}^k_{ml}$. This is justified from an asymptotic standpoint as all active variables should be found in models in ${\cal M}^k_{ml}$. Suppose there are $p_r$ variables remaining and $p_r< n$. Then, we can redefine the full model as the one with these $p_r$ variables and apply CMC to select a model from ${\cal M}^k_{ml}$. Simulation results show that when $n\gg k$, this approach works well.

\section{Appendix: proofs of lemmas}

\noindent {\bf Proof of Lemma \ref{lem1} }  For any fixed $\varepsilon >0$, there exists a $k$ such that $\alpha_k<\varepsilon$. By (A2), 
\[ \lim_{n\rightarrow \infty} P( \bs{\beta}^t \in  {\cal C}_{1-\alpha_k} ) = 1-\alpha_k> 1-\varepsilon. \]
It follows that there exists an $M>0$ such that for $n\geq M$, $P( \bs{\beta}^t \in  {\cal C}_{1-\alpha_k} ) > 1-\varepsilon$.
For $j>k$, since $\alpha_j<\alpha_k$ and thus $ {\cal C}_{1-\alpha_k} \subset  {\cal C}_{1-\alpha_j} $, we have 
\[P( \bs{\beta}^t \in  {\cal C}_{1-\alpha_j} ) \geq P( \bs{\beta}^t \in  {\cal C}_{1-\alpha_k} ) > 1-\varepsilon\] 
when $n\geq M$. 
Letting $N=\max\{k, M\}$, we have for $n\geq N$, $P( \bs{\beta}^t \in  {\cal C}_{1-\alpha_n} ) > 1-\varepsilon$.
This and $P( \bs{\beta}^t \in  {\cal C}_{1-\alpha_n} ) \leq 1$ imply (\ref{lemma1}).  \hfill $\Box$

\vspace{0.2in}

\noindent {\bf Proof of Lemma \ref{lem2}.}  For any $\bs{\beta} \in \mathbb{R}^{p+1}$, by triangle inequality,  
\beq
\| \bs{\beta} -\bs{\beta}^t\|_2 \leq \| \bs{\beta} -\hat{\bs{\beta}}\|_2 + \| \hat{\bs{\beta}}-\bs{\beta}^t \|_2.  \label{ineq}
\eeq
Under (A1) and (A3),  $\| \hat{\bs{\beta}}-\bs{\beta}^t \|_2=o_p(1)$ and $\max_{\bs{\beta}\in {\cal C}_{1-\alpha_n}} \| \bs{\beta} -\hat{\bs{\beta}}\|_2=o_p(1)$.  It follows from (\ref{ineq}) that
$\| \bs{\beta} -\bs{\beta}^t\|_2 = o_p(1)$ uniformly for all $\bs{\beta} \in {\cal C}_{1-\alpha_n}$ which implies (\ref{lemma2}). \hfill $\Box$

\vspace{0.5in}

\noindent {\bf \large Acknowledgement}

\vspace{0.1in}

This work is supported by the Natural Sciences and Engineering Research Council of Canada.
We thank Danielle Tsao of the University of Washington for helpful discussions and for assistance on a simulation study related to this paper.

\end{document}